\documentclass[11pt,a4paper]{article}
\usepackage[utf8]{inputenc}
\usepackage{amsmath,amssymb,amsthm}
\usepackage[margin=2.5cm]{geometry}
\usepackage{xcolor,booktabs,hyperref}
\usepackage{algorithm,algpseudocode}
\newtheorem{theorem}{Theorem}[section]
\newtheorem{proposition}[theorem]{Proposition}
\newtheorem{lemma}[theorem]{Lemma}
\newtheorem{corollary}[theorem]{Corollary}
\theoremstyle{definition}
\newtheorem{definition}[theorem]{Definition}
\newtheorem{remark}[theorem]{Remark}
\newtheorem{example}[theorem]{Example}
\newcommand{\X}{\mathcal{X}}
\newcommand{\K}{\Bbbk}
\newcommand{\TODO}[1]{\textcolor{red}{[TODO: #1]}}

\title{Singular rational $S$-functions of second-order ODEs:\\ structure, degree bounds and a Kovacic-type search}
\author{L.\,G.\,S.\ Duarte\thanks{Corresponding author. \texttt{lgsduarte@gmail.com}}
\and
L.\,A.\,C.\,P.\ da Mota\thanks{\texttt{lacpdamota@uerj.br}}\\[6pt]
\small Departamento de F\'{\i}sica Te\'orica, Instituto de F\'{\i}sica,\\
\small Universidade do Estado do Rio de Janeiro,\\
\small Rua S\~ao Francisco Xavier 524, 20550-013 Rio de Janeiro, RJ, Brazil}
\date{}

\begin{document}
\maketitle

\begin{abstract}
For a rational second-order ODE $y''=\phi(x,y,y')$, the $S$-function method reduces the search for a first
integral to finding a rational solution $S$ of the Riccati-type equation $\X(S)=S^2+\phi_zS-\phi_y$ along the
total-derivative field $\X$. The standard ansatz $S=P/N^k$, with $N=\operatorname{den}\phi$, misses every $S$
with poles outside $N=0$. Such \emph{singular} solutions cannot be avoided: for point transforms of linear
equations without hyperexponential solutions we describe every quasi-regular $S$ explicitly, and for a generic
transformation none exists. We then determine the structure of singular solutions. Outside $N$ every pole is
simple, with residue $-1$, and is not a Darboux surface; along non-characteristic factors of $N$ the leading
polar coefficient is free Cauchy data; and the only unknown denominator is a squarefree polynomial $D$, which
enters a linear equation. A Newton-polygon analysis at infinity gives the polynomial part up to a correction of
bounded weight. The weighted degree of $D$ is determined exactly when the excess is zero and the top Riccati
equation has no algebraic solution, and is bounded by residues of the top field otherwise, conditionally on the
vanishing of that correction. The dependence of $D$ on $x$ is obtained, without a degree bound, from the
rational solutions of a linear differential system. These results give a certified search procedure, implemented twice,
which finds the singular $S$ of nine non-Liouvillian pull-backs of the Airy
and Bessel equations. In one regime it is a decision procedure; for instance, it proves that $y''=xy+1$ has no
rational $S$.
\end{abstract}


(secondary).

\section{Introduction}

A rational second-order ODE $y''=\phi(x,y,y')$ that is not linearisable and has no Liouvillian first integral
lies beyond the reach of the classical Darboux and Prelle--Singer machinery. Typical examples are the pull-backs
of Airy's or Bessel's equation by a rational substitution, whose first integrals are built from Airy or Bessel
functions. The \emph{$S$-function method} of Duarte, da Mota and collaborators \cite{ACDM2019,NewS2023} reaches
part of this territory. It writes the first integral $I$ as $I_y=S\,I_z$ (with $z=y'$), where $S$ satisfies the
Riccati-type PDE
\[
\X(S)=S^2+\phi_zS-\phi_y,\qquad \X=\partial_x+z\partial_y+\phi\,\partial_z ,
\]
and, once a \emph{rational} $S$ is known, it reduces the problem to a first-order equation. A companion paper
\cite{PaperA} shows that, given a rational $S$, one can decide whether the integral is non-Liouvillian of
\emph{linear type}, and recover the underlying linear equation. The present paper is the front end of that
pipeline: it is about \emph{finding} the rational $S$.

In \cite{ACDM2019,NewS2023} the rational solutions are sought with the ansatz $S=P/N^k$, where $N$ is the
denominator of $\phi$; we call these solutions \emph{regular}. The definition in \cite[Def.~2.2]{NewS2023}
already admits denominators with factors outside $N$, but it does not say what those factors are, and the
ansatz cannot find them. We call such solutions \emph{singular}. Our first observation is that singular solutions
are not exotic. Non-Liouvillianity alone does not force them (Remark~\ref{rem:false}), but for point transforms of
linear equations without hyperexponential solutions every rational $S$ is singular unless the transformation has no
simple factor, and
in that exceptional case the regular $S$ is explicit and unique (Theorem~A, Section~\ref{sec:A}).

To find singular solutions we follow Kovacic's strategy for $\theta'+\theta^2=r$ \cite{Kovacic1986}: determine
the polar parts locally, determine the polynomial part at infinity, and recover the remaining global factor by
linear algebra once its degree is bounded. The transposition to~(\ref{ric}) is not mechanical. Three features of
the three-dimensional setting have no scalar analogue.
\begin{itemize}
\item \emph{Non-characteristic factors.} Along a factor $f$ of $N$ with $f_z\neq0$, the leading polar
coefficient of $S$ is free Cauchy data: no local analysis determines it (Proposition~\ref{prop:nonchar}). In one
variable this case does not occur, which is why Kovacic's algorithm is purely local. Here it forces the regular
ansatz $P/N_0$ into the method as a component, not an alternative.
\item \emph{The residue theorem kills the Fuchs count.} When the global factor $D$ is bounded through the
quasi-homogeneous top of the field at infinity, the analogue of Kovacic's exponent count is identically
satisfied, because the residues of $u/G$ sum to $1/g$. The degree of $D$ is bounded instead by the positivity
of local multiplicities. With the polynomial part known, the only obstruction is a polynomial first integral of
the top field, which plays the role of resonance (Theorem~B); in general the bound is conditional
(Remark~\ref{rem:Bscope}).
\item \emph{The top can itself be a Riccati equation.} On faces where all pieces of $\X$ have weight zero, the
top field is $\partial_x+G(x,u)\partial_u$ for a Riccati equation in $u=z/y$. When that equation is
non-Liouvillian, $-z/y$ is the unique weight-zero solution of the top problem, and the degree of $D$ comes out
exactly (Theorem~C). In this regime the non-Liouvillian case is the \emph{easy} one.
\end{itemize}

\paragraph{Contributions.}
\begin{enumerate}
\item Theorem~A: a complete description of the quasi-regular $S$ of point transforms of linear equations.
\item The local analysis: the residue lemma (Lemma~\ref{lem:res}), the characteristic/non-characteristic
dichotomy, and a denominator theorem (Theorem~\ref{thm:den}) that leaves a single unknown denominator $D$, which
enters linearly.
\item A structure theorem (Theorem~\ref{thm:struct}): every rational $S$ is a polynomial, plus a regular part
$P/N$, plus explicit terms $-\X(f)/f$ on non-characteristic factors, minus $\X(D)/D$.
\item The analysis at infinity by Newton polygons, with a decoupling of infinity from the finite poles up to a
correction of weight at most $e$ (Proposition~\ref{prop:decouple}).
\item Degree results for $D$: the weighted degree is exact when $e=0$ with a non-Liouvillian top (Theorem~C), and
bounded by residues when $e>0$, conditionally on the correction of Proposition~\ref{prop:decouple} vanishing
(Theorem~B).
\item A certified search procedure, Kovacic-$S$, implemented twice (Python/sympy and Maple). It has no false
positives; its completeness is proven only under the hypotheses of the theorems, and the remaining steps are
flagged as heuristic (Section~\ref{sec:algo}). In the regime of Theorem~C with $N_0=1$ it is a decision procedure
(Corollary~\ref{cor:Cfull}), which yields certified negative answers. We test it on non-Liouvillian pull-backs of
the Airy and Bessel equations.
\end{enumerate}

\paragraph{Scope.} We concentrate on the non-Liouvillian regime, where the regular ansatz fails and classical
Darboux methods do not apply. The complementary cases, a Liouvillian top or a resonant top, are exactly those in
which a Darboux structure is present. There the methods of \cite{ACDM2019,NewS2023} and the Darboux/Prelle--Singer
machinery already work; we discuss them in Section~\ref{sec:limits}.

\paragraph{Computer checks.} The statements rest on the proofs given here. In addition, the explicit identities
and all examples were checked by two independent implementations; scripts are in the supplementary material.

\section{Setting}\label{sec:setting}

\paragraph{The Riccati equation of $S$.}
Let $\K$ be a field of characteristic zero, and consider a rational second-order ODE
\begin{equation}\label{ode}
y''=\phi(x,y,y'),\qquad \phi=\frac{M}{N}\in\K(x,y,z),\quad \gcd(M,N)=1,
\end{equation}
where $z$ stands for $y'$. The total-derivative field is
$\X=\partial_x+z\,\partial_y+\phi\,\partial_z$, and a first integral of~(\ref{ode}) is a non-constant $I$ with
$\X(I)=0$. In the $S$-function approach \cite{ACDM2019,NewS2023}, one writes $I_y=S\,I_z$. Compatibility of the
mixed derivatives then forces $S$ to satisfy
\begin{equation}\label{ric}
\X(S)=S^2+\phi_zS-\phi_y ,
\end{equation}
and, conversely, a solution $S$ of~(\ref{ric}) reduces the problem to the first-order equation $dz/dy=-S$
(with $x$ as a parameter), the ``associated 1ODE'' of~\cite[Def.~1.4]{ACDM2019}. This paper is about finding the
\emph{rational} solutions of~(\ref{ric}).

\begin{definition}\label{def:regular}
A rational solution $S$ of~(\ref{ric}) is \emph{regular} if every irreducible factor of $\operatorname{den}S$
divides $N$. It is \emph{quasi-regular} if every irreducible factor of $\operatorname{den}S$ that depends on
$(y,z)$ divides $N$; factors in $x$ alone are allowed. Otherwise $S$ is \emph{singular}.
\end{definition}

The ansatz $S=P/N^k$ used in \cite{ACDM2019,NewS2023} finds exactly the regular solutions. The point of the
paper is that singular solutions are not exotic: for point transforms of linear equations without hyperexponential
solutions they are the rule (Section~\ref{sec:A}).

\paragraph{Linearisation.}
Writing $S=-\X(w)/w$ turns~(\ref{ric}) into the second-order linear equation
\begin{equation}\label{lin}
\X^2w-\phi_z\,\X w-\phi_y\,w=0 ,
\end{equation}
which is equation~(26) of~\cite{ACDM2019}. We use a relative version. Suppose $S=R-\X(D)/D$, where $R$ is a known
rational function and $D$ is an unknown polynomial. Then~(\ref{ric}) is equivalent to
\begin{equation}\label{LR}
\mathcal{L}_R(D):=\X^2D-(2R+\phi_z)\,\X D+\bigl(R^2-\X(R)+\phi_zR-\phi_y\bigr)D=0 ,
\end{equation}
which is \emph{linear} in $D$. The whole algorithm is organised so that, by the time $D$ is sought, $R$ is
already known up to finitely many choices. Then $D$ is found by linear algebra in a finite-dimensional space,
once its degree is bounded (Section~\ref{sec:bounds}).

\paragraph{A trivial branch.}
If $\phi_x=0$, then $S=-\phi/z$ always solves~(\ref{ric}): it corresponds to a first integral with $I_x=0$
\cite{ACDM2019}. In the language of~(\ref{LR}), $D=z$ always lies in the kernel. The algorithm flags this
branch and continues the search.

\section{Motivation: point transforms of linear equations}\label{sec:A}

In this section we determine all quasi-regular $S$ for the simplest class whose integrals can be
non-Liouvillian: point transforms, with $x$ fixed, of linear equations. The answer shows that, unless an
explicit algebraic condition holds, every $S$ is singular.

Let $L:\ y''=a(x)y'+b(x)y$ with $a,b\in\K(x)$, and let $T\in\K(x)[Y]$ with $n=\deg_YT\ge1$. The substitution
$y=T(x,Y)$, $z=T_x+T_YZ$ (with $Z=Y'$) transforms $L$ into
\begin{equation}\label{ET}
E_T:\quad Y''=\phi_T:=\frac{a(T_x+T_YZ)+bT-T_{xx}-2T_{xY}Z-T_{YY}Z^2}{T_Y},
\end{equation}
with field $\X_T=\partial_x+Z\partial_Y+\phi_T\partial_Z$. Write $g(x,Y,Z)=(x,T,T_x+T_YZ)$.

\begin{lemma}\label{lem:denT}
Every irreducible factor of $T_Y$ (as a polynomial in $Y$) divides $\operatorname{den}\phi_T$.
\end{lemma}
\begin{proof}
If $n=1$ there is nothing to prove. Otherwise, let $f^k\,\|\,T_Y$ with $k\ge1$. In characteristic zero,
$f^{k-1}\,\|\,T_{YY}$, so the coefficient $-T_{YY}/T_Y$ of $Z^2$ in $\phi_T$ has $f$ in its denominator with
exponent exactly one. Since $\phi_T$ is a polynomial in $Z$, this coefficient cannot be cancelled by the others.
\end{proof}

\begin{lemma}[transport]\label{lem:transport}
For arbitrary $\phi$ and $T$, set $\sigma:=S_T-\X_T(T_Y)/T_Y$. Then $S_T$ solves the Riccati equation of $E_T$
if and only if
\[
\X_T(\sigma)=\sigma^2+(\phi_z\circ g)\,\sigma-\phi_y\circ g .
\]
Moreover $\X_T(T)=T_x+T_YZ$ and $\X_T(T_x+T_YZ)=\phi\circ g$, that is, $\X_T$ and $\X$ are $g$-related.
\end{lemma}
\begin{proof}
For a rational $S$ put $\omega_S=-(\phi+zS)\,dx+S\,dy+dz$. A direct computation gives
$\omega_S\wedge d\omega_S=\bigl(\X(S)-S^2-\phi_zS+\phi_y\bigr)\,dx\wedge dy\wedge dz$, so $S$ solves~(\ref{ric}) if
and only if $\omega_S$ is integrable. Substituting $y=T$, $z=T_x+T_YZ$ and collecting terms,
\[
g^*\omega_{\sigma}=T_Y\,\omega^{T}_{S_T},\qquad S_T=\sigma\circ g+\frac{T_{xY}+T_{YY}Z}{T_Y},
\]
where $\omega^T$ is the corresponding form for $E_T$. Since $g^*(\omega\wedge d\omega)=g^*\omega\wedge d(g^*\omega)$,
multiplying a form by a nonzero function preserves integrability, and $g$ is dominant, $\omega_\sigma$ is integrable
if and only if $\omega^T_{S_T}$ is. The last two identities are immediate from $y=T$, $z=T_x+T_YZ$. (All identities
were also checked symbolically.)
\end{proof}

\begin{theorem}[Theorem A]\label{thm:A}
The quasi-regular solutions of the Riccati equation of $E_T$ are exactly the following.
\begin{enumerate}
\item[(i)] $S=s(x)+\X_T(T_Y)/T_Y$, where $s\in\K(x)$ satisfies $s'=s^2+as-b$. Equivalently, $v=e^{-\int s}$ is
a hyperexponential solution of $L$.
\item[(ii)] $S=-\dfrac{T_x+T_YZ+c'}{T+c}+\dfrac{\X_T(T_Y)}{T_Y}$, where $c\in\K(x)$ is a rational solution of $L$
($c=0$ allowed) and $T+c$ has \emph{no simple factor} in $Y$, that is, every irreducible factor of $T+c$ divides
$T_Y$.
\end{enumerate}
There is at most one $c\in\overline{\K}(x)$ for which $T+c$ has no simple factor, so branch~(ii) contributes at
most one $S$.
\end{theorem}

\begin{proof}
Let $\mathcal{R}=\K(x)[Y,1/T_Y]$. By Lemma~\ref{lem:denT}, $\operatorname{den}\phi_T$ divides $T_Y$ times a
polynomial in $x$. Hence a quasi-regular $S$ lies in $\mathcal{R}[Z]$, and so does $\sigma$, because
$\X_T(T_Y)=T_{xY}+T_{YY}Z$. For $L$, Lemma~\ref{lem:transport} reads $\X_T(\sigma)=\sigma^2+a\sigma-b$.

\emph{Degree in $Z$.} The field $\X_T$ maps $\mathcal{R}[Z]_{\le m}$ into $\mathcal{R}[Z]_{\le m+1}$. If
$m=\deg_Z\sigma\ge2$, the term $Z^{2m}$ occurs only in $\sigma^2$, which is impossible. Hence
$\sigma=p+q\,z$ with $z=T_x+T_YZ$ and $p,q\in\mathcal{R}$.

\emph{Coefficients in $z$.} Writing the equation in the independent coordinates $(x,Y,z)$ gives
\begin{align*}
z^2:&\quad q_Y=q^2T_Y,\\
z^1:&\quad p_Y/T_Y+q_x-q_YT_x/T_Y=2pq,\\
z^0:&\quad p_x-p_YT_x/T_Y=p^2+ap-b-bTq .
\end{align*}
The first equation gives $q=0$ or $q=-1/(T+c)$, with $c\in\K(x)$.
\begin{itemize}
\item If $q=0$, then $p_Y=0$ and $p'=p^2+ap-b$, which is case~(i).
\item If $q\neq0$, the $z^1$-equation becomes $((T+c)^2p)_Y=-c'T_Y$, so $p=-c'/(T+c)+e/(T+c)^2$ with
$e\in\K(x)$. Thus $\sigma=S_0\circ g$ for the \emph{rational} function
$S_0=-(z+c')/(y+c)+e/(y+c)^2$. Since $g$ is dominant, $S_0$ solves the Riccati equation of $L$. In powers of
$D=y+c$, the numerator of that equation is
$-L(c)\,D^3+(e'-ae)\,D^2-e^2$, where $L(c)=c''-ac'-bc$. Hence $e=0$ and $L(c)=0$.
\end{itemize}
Finally, $q=-1/(T+c)\in\mathcal{R}$ if and only if every irreducible factor of $T+c$ divides $T_Y=(T+c)_Y$, that
is, if and only if $T+c$ has no simple factor.

\emph{Converse.} The expressions in (i) and (ii) are the transports, by Lemma~\ref{lem:transport}, of the
solutions $s(x)$ and $-(z+c')/(y+c)$ of the Riccati equation of $L$. They are quasi-regular by
Lemma~\ref{lem:denT}.

\emph{Uniqueness of $c$.} If $c_1\ne c_2$, the roots of $T+c_1$ and of $T+c_2$ are disjoint. When neither
polynomial has a simple factor, each contributes total multiplicity at least $n/2$ to $T_Y$, which would force
$n\le n-1$.
\end{proof}

\begin{corollary}\label{cor:Tlin}
A linear equation $L$ has a quasi-regular $S$ if and only if it has a hyperexponential solution $v$. In that
case $S=-v'/v$.
\end{corollary}

\begin{corollary}\label{cor:nonhyper}
Suppose $L$ has no hyperexponential solution over $\overline{\K}(x)$ (so $L$ is not in Kovacic's case~1). Then $E_T$
has a quasi-regular $S$ if and only if $T$ has no simple factor in $Y$. In that case the quasi-regular $S$ is
unique,
\[
S_T=-\frac{T_x+T_YZ}{T}+\frac{\X_T(T_Y)}{T_Y}.
\]
In particular, if $T$ has a simple factor in $Y$ (the generic case), \emph{every} rational $S$ of $E_T$ is
singular, and singular solutions exist: the transport of $-z/y$, which solves the Riccati equation of $L$.
\end{corollary}

The condition ``$T$ has no simple factor'' is invariant under extension of $\K$, and it is tested by
$T/\gcd(T,T_Y)\mid T_Y$. Over a non-closed $\K$, what matters is the absence of $s\in\K(x)$. For instance, over
$\mathbb{Q}$ the equation $y''=-y$ has no $s\in\mathbb{Q}(x)$, so only branch~(ii) remains.

\begin{example}\label{ex:A}
Let $L$ be Airy's equation $y''=xy$ (or Bessel's equation of order zero), neither of which has a hyperexponential
solution.
\begin{itemize}
\item For $T=Y^3$, $Y^2(Y+1)^3$, $(Y^2+x)^2$ or $xY^2$, there is a regular $S_T$. For example, for $T=Y^3$ the
equation is $Y''=(xY^2-6Y'^2)/(3Y)$, with $S=-Z/Y$.
\item For $T=Y^2+x$, $Y^3+Y$ or $Y^3+xY+1$, every $S$ is singular. Airy with $T=Y^3+xY+1$ is the pull-back H4 of
the experiments (Section~\ref{sec:exp}). (H3 is \emph{not} a point transform of a linear equation: it is H1
composed with $Y=y^3+y$.)
\end{itemize}
\end{example}

\begin{remark}[two tempting but false statements]\label{rem:false}
``Non-Liouvillian first integral $\Rightarrow$ every $S$ is singular'' is false. The equation
$y''=x-(y'+y^2)^2-2yy'$ has the regular solution $S=2y$ ($N=1$), while $g=z+y^2$ satisfies $\X(g)=x-g^2$, so the
associated linear equation is Airy's and the integral is non-Liouvillian \cite{PaperA}. The weaker statement
``point transforms of linear equations without hyperexponential solutions have only singular $S$'' is also false
(take $T=Y^3$ in Example~\ref{ex:A}). Theorem~\ref{thm:A} says exactly when the singular factor of $-z/y$ is
absorbed into $N$: when $T$ has no simple factor.
\end{remark}

\begin{remark}[apparent zeros]\label{rem:apparent}
Regular and quasi-regular differ. The equation $y''=(3+(x-1)^2)y$ has the solution $v=(x-1)e^{x^2/2-x}$ and
$N=1$, and its only quasi-regular $S$ is $-v'/v=-1/(x-1)-(x-1)$, which is not regular. The factor $x-1$ is an
apparent zero of $v$, that is, a root of the polynomial $P$ in case~1 of Kovacic's algorithm. In the
terminology of~(\ref{LR}), it is the polynomial $D$: on linear equations, $D$ plays the role of Kovacic's
polynomial~$P$.
\end{remark}

\section{Local analysis along the factors of $N$}\label{sec:local}

Kovacic's algorithm for $\theta'+\theta^2=r$ is organised around one observation. At a pole of $\theta$ that is
not a pole of $r$, the pole is simple with residue~1, so all such poles together form $P'/P$ for a polynomial $P$;
at the poles of $r$ and at infinity, finitely many local candidates are computed from expansions. In this section
we transport the first half of that picture to~(\ref{ric}). The geometry of $(x,y,z)$-space adds a phenomenon
with no analogue in one variable: \emph{non-characteristic} factors of $N$, along which no local analysis
determines $S$. Throughout, $N=\prod f^{m_f}$ is the factorisation of $N$, and
$N_0:=\prod_{f_z\neq0}f^{m_f}$ is the part of $N$ made of non-characteristic factors (Section~\ref{sec:nonchar}).

\subsection{Poles outside $N$}

\begin{lemma}[residue lemma]\label{lem:res}
Let $S=P/Q$, in lowest terms, be a rational solution of~(\ref{ric}), and let $f$ be an irreducible factor of $Q$
with $f\nmid N$. Then:
\begin{enumerate}
\item[(a)] $f$ has multiplicity one in $Q$;
\item[(b)] $f$ is \emph{not} a Darboux polynomial of $\X$;
\item[(c)] $S+\X(f)/f$ has no pole along $f=0$; that is, the residue of $S$ along $f$ is $-1$.
\end{enumerate}
\end{lemma}
\begin{proof}
Along $f=0$ the coefficients of $\X$ are regular, since $f\nmid N$. Let $k\ge1$ be the order of the pole, with
$S=P_0/f^k+\cdots$ and $P_0\not\equiv0\pmod f$. Then:
\begin{itemize}
\item $\X(S)$ has a pole of order at most $k+1$, with principal part $-kP_0\X(f)/f^{k+1}$;
\item $S^2$ has a pole of order exactly $2k$;
\item $\phi_zS$ and $\phi_y$ have poles of order at most $k$.
\end{itemize}
If $k\ge2$, the principal part of $S^2$ cannot be cancelled, which proves~(a). For $k=1$, the order-two terms
give $-P_0\X(f)\equiv P_0^2\pmod f$, so $P_0\equiv-\X(f)\pmod f$, which is~(c). If $f$ were Darboux,
$\X(f)\equiv0\pmod f$ would force $P_0^2\equiv0$, a contradiction; this is~(b).
\end{proof}

\begin{remark}
Parts~(a) and~(c) are the classical scalar fact, transported along $\X$: a pole of a rational solution of a
Riccati equation that is not a pole of the coefficients is simple, with residue fixed by the normalisation
(Kovacic's analysis of case~1, \cite[Sec.~3]{Kovacic1986}; see also \cite{Echi}). Part~(b)
has no scalar analogue: it says that the unknown polynomial collecting these poles is never an invariant surface
of the flow. Consequently, the poles of $S$ outside $N$ are collected in $-\X(D)/D$ for a squarefree polynomial
$D$, coprime to $N$, with no Darboux factor.
\end{remark}

\subsection{Factors of $N$: the characteristic case $f_z=0$}

Let $f$ be an irreducible factor of $N$ of multiplicity $m$. Write $N=f^mN'$, $\mu=M/N'$ (regular and nonzero
along $f=0$) and $X_0=\partial_x+z\partial_y$. Let $k\ge0$ be the order of the pole of $S$ along $f$, with
$S=s/f^k+\cdots$ and $s$ taken modulo $f$. We assume $f_y\not\equiv0\pmod f$; factors depending on $x$ alone are
treated separately (Remark~\ref{rem:xonly}).

\begin{proposition}\label{prop:char}
If $f_z=0$, then $k\le m$, and modulo $f$:
\begin{align*}
&\text{(E1)}\quad k=m=1: && \mu\,s_z=s^2+\bigl(\mu_z+X_0(f)\bigr)s+\mu f_y,\\
&\text{(E2)}\quad k=m\ge2: && \mu\,s_z-\mu_zs=s^2,\quad\text{i.e. } s=\mu/(c-z),\\
&\text{(E3)}\quad k=1<m: && \mu\,s_z-\mu_zs=m\,\mu f_y,\quad\text{i.e. } s=\mu\bigl(c+mf_y\textstyle\int dz/\mu\bigr),\\
&\text{(E3$'$)}\quad 2\le k<m: && \mu\,s_z-\mu_zs=0,\quad\text{i.e. } s=c\,\mu,
\end{align*}
where $c\in\K(x,y)$ is taken modulo $f$.
\end{proposition}
\begin{proof}
Since $f_z=0$, $X_0(f)$ is regular. The pole orders along $f$ are: $k+1$ for $-ks\,X_0(f)/f^{k+1}$; $m+k$ for
$\phi S_z$ and $\phi_zS$ (leading coefficients $\mu s_z$ and $\mu_zs$); $2k$ for $S^2$; and $m+1$ for $\phi_y$
(leading coefficient $-m\mu f_y$). If $k>m$, then $S^2$ is alone at the top order. Equating the top orders in the
remaining cases gives (E1)--(E3$'$). In (E2), $\sigma=1/s$ satisfies $-(\mu\sigma)_z=1$.
\end{proof}

Equation (E1) is a Riccati equation in the \emph{single} variable $z$, with coefficients in the function field of
the curve $f=0$: the problem along $N$ drops one dimension. Suppose the polar part of $S$ along $f$ is logarithmic,
$\kappa\,\X(f)/f$, and \emph{the global factor does not contribute to it}, that is,
\begin{equation}\label{Hlog}
\text{(H$_{\log}$)}\qquad
\begin{minipage}{0.8\linewidth}
along every characteristic factor $f\mid N$ (including factors in $x$ alone), either $S$ has no pole, or
$m_f=1$, the polar part of $S$ is $\kappa_f\X(f)/f$ with $\kappa_f$ a rational constant, and $f\mid D_z$.
\end{minipage}
\end{equation}
Then $s=\kappa X_0(f)$, and substituting into (E1) with $X_0(f)_z=f_y$ gives the \emph{indicial equation}
\begin{equation}\label{indicial}
(\kappa+1)\kappa\,X_0(f)^2+\kappa\,\mu_z\,X_0(f)-(\kappa-1)\,\mu f_y\equiv0\pmod f .
\end{equation}
This is polynomial in $z$ over the curve, and its coefficients are equations for $\kappa$. It is the analogue of
Kovacic's local exponents. On battery H, where every $D$ is free of $z$, it returns the true exponents (for instance,
$\kappa=1$ on the factor $x+3y^2$ of H4, and $\kappa\in\{-2,-1\}$ on the factor $y+1$ of H5, the true one being
$-1$).

\begin{remark}[when (H$_{\log}$) fails]\label{rem:Hlog}
In general $-\X(D)/D$ has a pole on $f$ through $\phi D_z/D$, and the polar coefficient is
$s\equiv\kappa X_0(f)-\mu D_z/D\pmod f$. Then $\kappa$ is not determined locally. For example, $y''=y'^2(y'+2)/y$ has
the solution $S=\X(y)/y-\X(xz+y)/(xz+y)$, with $\kappa=1$ on $f=y$ and $D=xz+y$, but~(\ref{indicial}) has no root.
Similarly, characteristic factors with $m_f\ge2$ can carry poles of order up to $m_f$: for
$y''=(2xy'+yy'^2-y)/y^2$ one has $S=-(2x+yz)/y^2$, which is case (E2) with $c=0$. Both examples are Liouvillian, but
nothing in this section uses that. A complete local step must treat the coupled condition, or
solve (E1)--(E3$'$) with a rational-Riccati algorithm over the function field of the curve; we leave this open and
state (H$_{\log}$) as a hypothesis wherever the indicial equation is used. All these solutions still have the form of
Theorem~\ref{thm:struct}, so the structural step of the search covers them. The options it excludes are: characteristic
factors with $m_f\ge2$ and poles of order up to $m_f$; non-logarithmic or non-constant polar coefficients along
characteristic factors (as in $\phi=-(2xyz^2-y^2z+2yz^3-1)/(xy^2-1)$, $S=2yz(x+z)/(xy^2-1)$, which has $D=1$); and
the polar part $-\X(f)/f$ on non-characteristic factors (Proposition~\ref{prop:nonchar}(d)).
\end{remark}

\subsection{Factors of $N$: the non-characteristic case $f_z\neq0$}\label{sec:nonchar}

Near a generic point of $f=0$, use coordinates $(x,y,w)$ with $w=f$ and $z=Z(x,y,w)$. Put $C=\mu f_z$ and
expand $S=\sum_{j\ge0}s_j(x,y)\,w^{j-k}$.

\begin{proposition}\label{prop:nonchar}
Let $f\mid N$ be irreducible with $f_z\neq0$ and multiplicity $m$. Then:
\begin{enumerate}
\item[(a)] $C\not\equiv0\pmod f$;
\item[(b)] $1\le k<m$ is impossible;
\item[(c)] if $k=m$, the equation imposes \emph{nothing} on $s_0$; each $s_j$ with $j\ge1$ is determined by
$s_0,\dots,s_{j-1}$, with coefficient $jC_0$;
\item[(d)] if $k=m+1$, then $s_0=-C_0$, and the whole polar part is determined. This is the principal part of
$-\X(f)/f$;
\item[(e)] $k>m+1$ is impossible.
\end{enumerate}
\end{proposition}
\begin{proof}
(a) $f$ is irreducible and $\deg_zf_z<\deg_zf$, so $f\nmid f_z$. Also $\mu N'=M$ and $\gcd(M,N)=1$, so
$f\nmid\mu$.

In the coordinates $(x,y,w)$, $\X=\partial_x+Z\partial_y+\bigl(X_0(f)+C/w^m\bigr)\partial_w$. Hence $\X(f)$ has a
pole of order $m$ with leading coefficient $C_0$, and $\phi_z$ has a pole of order $m+1$ with leading coefficient
$-mC_0$.

(b) At order $m+k+1$, only $-kC_0s_0$ (from $\X(S)$) and $-mC_0s_0$ (from $\phi_zS$) appear, plus $s_0^2$ when
$2k=m+k+1$. For $1\le k<m$ this gives $(m-k)C_0s_0\equiv0$, which is impossible by~(a).

(c), (d) The top order vanishes identically when $k=m$. At order $2m+1-j$, the coefficient $s_j$ enters linearly,
with coefficient $jC_0$. For $k=m+1$, the top order gives $s_0(s_0+C_0)=0$, and the subsequent $s_j$ enter with
coefficient $(j+1)C_0\neq0$.

(e) If $k>m+1$, then $2k>m+k+1$, so $s_0^2$ stands alone at the top order $2k$.
\end{proof}

The two cases deserve names. A factor with $f_z=0$ is \emph{characteristic}: the surface $f=0$ is a union of
trajectories in the $z$-direction, and local analysis yields finitely many exponents~(\ref{indicial}). A factor
with $f_z\ne0$ is \emph{non-characteristic}. Proposition~\ref{prop:nonchar}(c) says that the leading coefficient
$s_0$ is free Cauchy data on the surface: \emph{no local analysis can determine it}. In one variable the second
case does not exist, which is why the classical algorithm is purely local. In three variables, hypersurfaces
transversal to $\partial_z$ along which $\phi$ blows up release a whole function's worth of freedom, and it must
be recovered by a global computation (Section~\ref{sec:infinity}).

\begin{remark}
With $N_0=f$ irreducible and $m=1$, the Cauchy data is $s_0\equiv P_0-\mu D_z/D\pmod f$ in the ansatz of
Section~\ref{sec:infinity}. Its only pole on $f=0$ lies on the curve $\{f=0\}\cap\{D=0\}$, with order one. This
identity is the ansatz rewritten on the surface, not a new condition, but it describes the search space exactly
(checked on the cases C5 and C6 of the blind battery of Section~\ref{sec:robust}, where $s_0=-2z/(y+z)$ and
$-4z/(y+z)$, both with $D=y+z$).
\end{remark}

\begin{remark}\label{rem:xonly}
Factors $f=f(x)$ of $N$, with $f_y=f_z=0$, are excluded from Proposition~\ref{prop:char}; they occur in the
experiments (H6 has $N=x^2$, up to a constant). The same order count, now with $X_0(f)=f'(x)$ regular and $\phi_y$ of order only $m$,
gives $k\le m$ and, modulo $f$,
\[
k=m=1:\ \ \mu\,s_z=s^2+\bigl(\mu_z+f'\bigr)s;\qquad k=m\ge2:\ \ \mu\,s_z-\mu_zs=s^2;\qquad 1\le k<m:\ \ s=c\,\mu ,
\]
with $c\in\K(x,y)$ taken modulo $f$ (checked symbolically for $(m,k)=(1,1),(2,2),(2,1),(3,2)$). There is no
$f_y$-term. For a logarithmic polar part $\kappa\X(f)/f=\kappa f'/f$ with $k=m=1$, we have $s=\kappa f'$ and $s_z=0$,
so the first equation becomes the indicial condition $\kappa(\kappa+1)f'^2+\kappa\,\mu_zf'\equiv0\pmod f$.
Since it must hold identically in $z$ and $f\nmid f'$, for $\kappa\neq0$ it forces $\mu_z\equiv-(\kappa+1)f'\pmod f$:
$\mu$ is affine in $z$ on $f=0$. The implementation
allows such factors in the coefficients of $D$ and $R$ rather than treating them by this analysis.
\end{remark}

\subsection{The denominator theorem}

\begin{theorem}\label{thm:den}
Let $S$ be a rational solution of~(\ref{ric}) and $N=\prod f^{m_f}$. Then
\[
\operatorname{den}S\ \Big|\ D\cdot N\cdot\prod_{f\mid N,\ f_z\ne0}f ,
\]
where $D$ is squarefree, coprime to $N$, has no Darboux factor, and $S+\X(D)/D$ has no pole on $D=0$. On a
non-characteristic factor, the extra order occurs only through the term $-\X(f)/f$.
\end{theorem}
\begin{proof}
Combine Lemma~\ref{lem:res} with Proposition~\ref{prop:char} and Remark~\ref{rem:xonly} ($k\le m$ when
$f_z=0$) and Proposition~\ref{prop:nonchar}(b),(e) ($k\in\{0,m,m+1\}$ when $f_z\neq0$).
\end{proof}

In words, the only unknown denominator is $D$, and it enters~(\ref{LR}) linearly. Everything else in the
denominator of $S$ is read off from $N$. This answers the question left open by the definition
in~\cite[Def.~2.2]{NewS2023}: the framework there allows factors beyond $N$ (one checks that the denominator
used there divides $N$) but does not say what they are.

\section{Analysis at infinity and the decoupling}\label{sec:infinity}

\subsection{Weights and dominant balance}
For a weight vector $w=(w_x,w_y,w_z)$, the three pieces of $\X$ have weights
\[
w(\partial_x)=-w_x,\qquad w(z\partial_y)=w_z-w_y,\qquad w(\phi\,\partial_z)=w(\phi)-w_z ,
\]
where $w(\phi)=w(M)-w(N)$. Let $e$ be the largest of the three, and let $\X_{\rm top}$ be the sum of the pieces
that attain it. If $S$ has weight $\sigma$, the four terms of~(\ref{ric}) sit at levels
$e+\sigma$, $2\sigma$, $(w(\phi)-w_z)+\sigma$ and $w(\phi)-w_y$. A dominant balance in which $S^2$ takes part
requires
\begin{equation}\label{balance}
\sigma\in\Bigl\{\,e,\ \tfrac12\bigl(w(\phi)-w_y\bigr),\ w(\phi)-w_y-e\,\Bigr\}.
\end{equation}
There is one more possibility, the \emph{degenerate top}: the dominant level contains only the linear terms
$\X_{\rm top}(S_{\rm top})$ and $(\phi_z)_{\rm top}S_{\rm top}$ (this needs $\sigma<e$, and
$w(\phi)-w_y<e+\sigma$). Then $\sigma$ is not fixed by~(\ref{balance}) but by the weights of the solutions of a
linear homogeneous equation. On a balanced face (defined below) all three values in~(\ref{balance}) equal $e$, and
$\phi_y$ sits at level $2e$; so a top with $\sigma<e$ is possible only if $(\phi_y)_{\rm top}$ vanishes. With
$\phi_{\rm top}=y^{2\omega-1}F(u)$, $u=z/y^\omega$ (Section~\ref{sec:bounds}), this means $(2\omega-1)F=\omega uF'$,
whose only solutions are $F=C\,u^{(2\omega-1)/\omega}$, not rational when $\omega>1$. Hence, on a balanced face with
$e>0$, \emph{every} rational solution of~(\ref{ric}) has weight exactly $e$: for $\sigma>e$ the term $S^2$ would be
alone at the top, and for $\sigma<e$ the term $\phi_y$.

The relevant weights are not taken from a grid. They come from the Newton polytope of $(M,N)$: with $w_x,w_y$ fixed,
$w(\phi)$ is a piecewise-affine function of $w_z$, and the candidates for $w_z$ are its breakpoints (where two
monomials of $M$ or of $N$ tie at the top) together with the solutions of $w(\phi)=2w_z-w_y$ on each affine piece.
This gives a finite list. The \emph{balanced} weights, $w(\phi)=2w_z-w_y$, where $z\partial_y$ and $\phi\partial_z$
have the same weight, are tried first; they are the ones relevant to Section~\ref{sec:bounds}. (For $y''=y^3$, for
instance, $\phi$ is a single monomial and the balanced weight $w_z=2$ is not a breakpoint.)

\subsection{The top equation and the linear recursion}
For each admissible $\sigma$, the dominant level of~(\ref{ric}) is an equation for the top part $S_{\rm top}$.
We solve it with an ansatz of quasi-homogeneous Laurent monomials of weight $\sigma$; the system is small and
quadratic. This ansatz is heuristic: the top problem is a Riccati-type equation (in $u=z/y^\omega$, see
Section~\ref{sec:bounds}), whose rational solutions need not be monomials. Below the top, the expansion is computed
weight by weight, and each step is \emph{linear}: the new term $t$ satisfies $\mathcal{T}(t)=(\text{known})$ with
$\mathcal{T}=\X_{\rm top}-2S_{\rm top}-(\phi_z)_{\rm top}$. The polynomial part of the result is denoted $T$. If
$\mathcal{T}$ has a kernel at some weight (a resonance, as in Painlev\'e--Kovalevskaya analysis), $T$ is determined
only up to free parameters; the implementation does not treat this case, which is one more heuristic step. For a
degenerate top the top equation is linear and homogeneous, so $S_{\rm top}$ is determined only up to a factor, fixed
at the first level where $S^2$ enters.

\begin{example}
For the pull-back of $y''=xy-y^5-4y^2y'$ by $y=Y^3+xY+1$, the face $(w_y,w_z)=(1,7)$ gives
$S_{\rm top}=2y^6-z/y$, uniquely. The linear recursion then reproduces, weight by weight from 5 down to 0,
\[
T=2y^6+4xy^4+4y^3+2x^2y^2+4xy+2=2(y^3+xy+1)^2 ,
\]
which is exactly the polynomial part of the true $S$, obtained \emph{without knowing $D$}.
\end{example}

\subsection{The ansatz and the decoupling}
\begin{theorem}[structure]\label{thm:struct}
Every rational solution of~(\ref{ric}) can be written as
\begin{equation}\label{struct}
S=Q+\frac{P}{N}-\sum_{f\in\mathcal{E}}\frac{\X(f)}{f}-\frac{\X(D)}{D},
\end{equation}
where:
\begin{itemize}
\item $Q\in\K(x)[y,z]$ is a polynomial in $(y,z)$;
\item $P\in\K(x)[y,z]$ (the normalisation of $P/N$ relative to a weight is fixed after Proposition~\ref{prop:decouple});
\item $\mathcal{E}$ is a set of non-characteristic factors of $N$, namely those along which $S$ has a pole of order $m_f+1$;
\item $D$ is squarefree, coprime to $N$, with no Darboux factor, and collects the poles of $S$ outside $N$.
\end{itemize}
\end{theorem}
\begin{proof}
By Lemma~\ref{lem:res}, $S+\X(D)/D$ has no pole outside $N$. By Proposition~\ref{prop:char} (and
Remark~\ref{rem:xonly}), along a characteristic factor the pole order is at most $m_f$. By
Proposition~\ref{prop:nonchar}, along a non-characteristic factor it is $0$, $m_f$ or $m_f+1$. In the last case the
leading coefficient is $-C_0$, which is the leading coefficient of $-\X(f)/f$, so $S+\X(f)/f$ has order at most
$m_f$ there. What remains has denominator dividing $N$ times a polynomial in $x$; write it as $Q+P/N$ (for instance with $Q=0$).
\end{proof}
The theorem holds without (H$_{\log}$), and we checked it on every example of this paper, including those of
Remark~\ref{rem:Hlog} and the 14 cases of Section~\ref{sec:coverage}. The logarithmic terms
$\kappa_f\X(f)/f=\kappa_fX_0(f)/f$ are the special case of $P/N$ used in the linear step. The term $P/N$ is exactly the regular ansatz of~\cite{ACDM2019,NewS2023}, which is therefore
not an alternative to the present method but a part of it: there are singular solutions that need both $P/N$ and
$D$. The search therefore has two kinds of steps:
\begin{itemize}
\item the \emph{linear step}, under (H$_{\log}$) and $N_0=1$: the regular part is $\sum\kappa_f\X(f)/f$ with
$\kappa_f$ from the indicial equation, $R$ is formed as in Algorithm~\ref{alg:KS}, the equation for $D$ is linear,
and Theorems~\ref{thm:B} and~\ref{thm:C} apply;
\item the \emph{structural step}: $P$ general (with the normal form below), which is bilinear in $(P,D)$.
\end{itemize}
An example of a singular solution needing both parts, case C5 of Section~\ref{sec:robust}, is
\begin{equation}\label{C5}
\phi=\frac{-y^2z-y-z^3}{y^2+z^2+1},\qquad S=1-\frac{1}{N}-\frac{\X(y+z)}{y+z},\qquad N=y^2+z^2+1 .
\end{equation}
The system~(\ref{struct}) is linear in $D$ and quadratic in $P$. The next observation removes most of the
unknowns of $P$.

Fix a weight $w$. Every rational function $F$ has a unique expansion at infinity
$F=\sum_{k\le k_0}F_{[k]}$ into quasi-homogeneous pieces. With $u=z/y^{\omega}$, each piece is $y^k r_k(u)$, and
splitting $r_k$ into its polynomial and proper parts in $u$ decomposes every piece uniquely into Laurent monomials
and a proper part. We write $\operatorname{pp}_w(F)$ for the sum of the Laurent monomials with non-negative
exponents in $(y,z)$; this is the \emph{polynomial part of $F$ at infinity}.

\begin{proposition}[decoupling]\label{prop:decouple}
Assume $N_0=1$, and write $S=Q+\sum_f\kappa_f\X(f)/f-\X(D)/D$ with $Q$ polynomial. Let $T=\operatorname{pp}_w(S)$,
which is computed by the Newton recursion. Then
\[
Q=T-\operatorname{pp}_w\Bigl(\sum_f\kappa_f\frac{\X(f)}{f}\Bigr)+\operatorname{pp}_w\Bigl(\frac{\X(D)}{D}\Bigr),
\]
where the middle term is known, and the last term $c:=\operatorname{pp}_w(\X(D)/D)$ is a polynomial whose weights
lie in $[0,e]$. In particular, every part of $Q$ of weight greater than $e$ is known.
\end{proposition}
\begin{proof}
$\operatorname{pp}_w$ is linear and is the identity on polynomials. Since $\X$ raises weights by at most $e$,
$\X(D)/D$ has weight at most $e$.
\end{proof}

\begin{example}[the correction $c$ is needed]\label{ex:gap}
Let $\phi=-y^5-4y^2z$ (face $\omega=3$, $e=2$) and $D=2xz+y$. Then $S=-\X(D)/D$ solves~(\ref{ric}), with $Q=0$ and no
$\kappa$. However,
\[
S=\frac{2xy^5+8xy^2z-3z}{2xz+y},\qquad \operatorname{pp}_w(S)=4y^2-\frac{3}{2x}\neq0 .
\]
With $R=T$, the linear equation $\mathcal{L}_T(D)=0$ has \emph{no} nonzero solution of weight $\le3$ (checked with
coefficients of degree $\le1$ in $x$). With the
correction $c=c_0(x)+c_1(x)y+c_2(x)y^2$ left unknown, the bilinear system in $(c,D)$ returns $c=-T$ (so $R=0$) and
the family $D=2xz+y+\lambda z$, which contains the trivial branch (a second component, $D=xz$ with $R=1/x$,
gives $S=-\phi/z$ again). This took 2.3\,s in sympy and 0.06\,s in Maple, with the same result.
\end{example}

Since $D$ is only defined up to a factor in $\K(x)$, a weight-zero term $h_0'/h_0$ of $c$ coming from
$D_{\rm top}=y^dh_0(x)$ can be absorbed into $D$. On battery~H, $c=0$ holds a posteriori, which explains why the
linear step sufficed there. In general the algorithm runs the linear step with $c=0$ and then a \emph{correction
step} with $c$ unknown, which is bilinear in $(c,D)$ but involves only the few monomials of weight in $[0,e]$
(Algorithm~\ref{alg:KS}).

\begin{remark}[on balanced faces the top ansatz loses nothing]\label{rem:balanced}
On a balanced face with $e>0$, $S$ has weight exactly $e$ (Section~\ref{sec:infinity}), so $T=\operatorname{pp}_w(S)$,
$\operatorname{pp}_w(\sum\kappa_f\X(f)/f)$ and the correction $c$ all have weights in $[0,e]$. The whole polynomial
part $Q$ of~(\ref{struct}) is then a polynomial of weight at most $e$, which the correction step treats as unknown
when it is run from the candidate $T=0$ (always included). So a top that the Laurent-monomial ansatz misses cannot
make the search lose a solution on such a face; it only moves the solution from the linear step to the bilinear
correction step, whose degree limit is heuristic (L4). The ansatz matters for speed, and for the bound of
Theorem~\ref{thm:B}, which uses $R_{[e]}$.

This happens. By Theorem~\ref{thm:B}(1), $D_{\rm top}=y^dh(u)$ where $h$ solves a linear ODE in $u$, and $h$ may
vanish at an \emph{ordinary} point $u=a$ ($G(a)\neq0$); then $S_{\rm top}$ has a pole at $u=a$ and is not a sum of
Laurent monomials. (Zeros at the roots $u_i$ of $G$ are harmless, since $G/(u-u_i)$ is a polynomial.) Imposing
$h=u-a$ on the top equation~(\ref{equ}) with $\omega=3$ and $R_{[e]}=t\,y^2$ gives, among others, $F=-8a^2+10au$,
$t=2a$, $d=4$, with $G(a)=-a^2\neq0$. For $a=1$ the full equation is
\[
y''=-8y^5+10y^2y',\qquad S=2y^2-\frac{\X(D)}{D},\quad D=y\,(z-y^3),\qquad
S_{\rm top}=-y^2\,\frac{u^2+4u-6}{u-1}.
\]
Here the ansatz produces no top with a nonzero polynomial part, and the linear step with $T=0$ returns the kernel
spanned by $z$ (the trivial branch) and $2xz+y$, and no member of the resulting family is $S$. The correction step with
$Q=q_0+q_1y+q_2y^2$ and $w\deg D\le4$ returns exactly $Q=2y^2$, $D=y(z-y^3)$, as its only component, in
$0.15$\,s (sympy) and $0.03$\,s (Maple~17); both implementations also return the same kernel for the linear
step. Note that
$\operatorname{pp}_w(S)=-5y^2$, so the correction is $c=7y^2\neq0$. The equation is autonomous and
quasi-homogeneous, hence Liouvillian, but the phenomenon lives entirely at the top level.
\end{remark}

\begin{remark}[normalisation and $N_0\neq1$]\label{rem:N0}
In Theorem~\ref{thm:struct}, $P/N$ can always be normalised by $\operatorname{pp}_w(P/N)=0$, moving
$\operatorname{pp}_w(P/N)$ into $Q$. With this normalisation the proof of Proposition~\ref{prop:decouple} applies
verbatim to an arbitrary $N$:
\[
Q=T+\operatorname{pp}_w\Bigl(\sum_{f\in\mathcal{E}}\frac{\X(f)}{f}\Bigr)+c ,\qquad
c=\operatorname{pp}_w\Bigl(\frac{\X(D)}{D}\Bigr),
\]
so again every part of $Q$ of weight greater than $e$ is known, and the unknowns are $P$ (subject to the linear
condition $\operatorname{pp}_w(P/N)=0$), the correction $c$ and $D$. No degree bound for this case is proved here.
\end{remark}

The content of the proposition is Kovacic's philosophy: \emph{the analysis at infinity and the analysis at finite
poles decouple up to a correction of weight at most $e$}, and all local data, including the data at infinity, are
exhausted before anything global is solved. Its effect is large. For~(\ref{C5}) the unknowns of $P$ drop from 9 to
6, and the solving time drops from more than 40 minutes to about 2--3 seconds (sympy; Table~\ref{tab:bilinear}).
With Maple's \texttt{Groebner[Basis]}, the reduced bilinear step takes under 0.1\,s for $\deg D\le2$ and under
8\,s on the hardest cases tested (Section~\ref{sec:exp}). When $N_0=1$ and $c=0$, which is the case throughout
battery~H, the step is purely linear.

\section{Degree bounds for $D$}\label{sec:bounds}

By~(\ref{LR}), once $R=S+\X(D)/D$ is fixed (in the linear step,
$R=T-\operatorname{pp}_w(\sum\kappa_f\X(f)/f)+\sum\kappa_f\X(f)/f$), $D$ ranges over the kernel of a linear operator,
and only a bound on $\deg D$ is missing. As in Kovacic's algorithm, the bound comes from infinity. Two regimes
appear, according to the excess $e$ of a balanced face. In both we assume $N_0=1$ and (H$_{\log}$). The theorems are stated for the true $R$; how much of $R$ the algorithm
knows in advance is discussed in Remark~\ref{rem:Bscope}.

\subsection{The regime $e>0$: a bound by residues}
Let $w=(0,1,\omega)$ be balanced, $w(\phi)=2\omega-1$, with $e=\omega-1>0$. Assume that $N_{\rm top}=\nu(x)y^\ell$
is a monomial in $(y,z)$. With $u=z/y^\omega$, write $\phi_{\rm top}=y^{2\omega-1}F(u)$, and assume $F$
is polynomial of degree at most $2$ with $g:=f_2-\omega\ne0$, where $f_2$ is the coefficient of $u^2$. Then
\[
\X_{\rm top}=y^e\bigl[u\,y\,\partial_y+G(u)\,\partial_u\bigr],\qquad G=F-\omega u^2\sim g\,u^2 .
\]
Write $D_{\rm top}=y^dh(u)$, $\theta_dh=d\,u\,h+Gh'$, $(2R+\phi_z)_{[e]}=y^e\alpha(u)$ and
$(R^2-\X(R)+\phi_zR-\phi_y)_{[2e]}=y^{2e}\beta(u)$, and assume $\deg_u\alpha\le1$, $\deg_u\beta\le2$, with
$\alpha_\infty,\beta_\infty$ the coefficients of $u$ and $u^2$. Then the top level of~(\ref{LR}) is the ODE in $u$
\begin{equation}\label{equ}
\theta_{d+e}\theta_dh-\alpha\,\theta_dh+\beta\,h=0 .
\end{equation}

\begin{theorem}[Theorem B]\label{thm:B}
Assume moreover that $G$ has simple roots $u_i$ and that $\alpha,\beta$ are regular at them. Let
$r_i=\operatorname{Res}_{u_i}(u/G)=u_i/G'(u_i)$, let $\rho_{i,\pm}$ be the roots of
$\rho(\rho+eu_i)-\alpha(u_i)\rho+\beta(u_i)=0$, and put $c_{i,\pm}=\rho_{i,\pm}/G'(u_i)$. Call $d$
\emph{admissible} if~(\ref{equ}) has a nonzero polynomial solution $y^dh$. Then, for admissible $d$:
\begin{enumerate}
\item the multiplicity of $z-u_iy^\omega$ in $D_{\rm top}$ is $n_i=c_{i,\pm}-d\,r_i$ for one choice of sign;
\item $\sum_ir_i=1/g$;
\item if some $r_i$ is a positive constant, then $d\le\bigl\lfloor\min_{r_i>0}\max_\pm c_{i,\pm}/r_i\bigr\rfloor$,
the maximum being taken over real constant branches;
\item if $1+\omega/g<0$, then $d\le\omega\,\sigma_{\max}/(\omega+g)$, where $\sigma$ ranges over the
roots of $(\sigma+e+g)\sigma-\alpha_\infty\sigma+\beta_\infty=0$;
\item if some $r_i\notin\mathbb{Q}$ (or is not constant), there is at most one $d$ for each choice of signs;
\item the set of admissible $d$ is infinite if and only if $\X_{\rm top}$ has a non-constant polynomial first
integral, namely a power of $H=y^{1+\omega/g}\prod_i(z-u_iy^\omega)^{-r_i}$, and at least one $d$ is admissible.
\end{enumerate}
Since $w_y=1<\omega$, the total degree of $D$ in $(y,z)$ is at most $d$.
\end{theorem}

\begin{proof}
(1) Put $h=(u-u_i)^n\tilde h$ with $\tilde h(u_i)\ne0$. Then $\theta_dh=\lambda(u-u_i)^n\tilde h+\cdots$ and
$\theta_{d+e}\theta_dh=(\lambda+eu_i)\lambda(u-u_i)^n\tilde h+\cdots$, with $\lambda=d\,u_i+nG'(u_i)$. The lowest
order of~(\ref{equ}) gives $\lambda=\rho_{i,\pm}$.

(2) The poles of $u/G$ are the $u_i$ and $\infty$, and $\operatorname{Res}_\infty(u/G)=-1/g$.

(3) and (5) follow from $n_i\in\mathbb{Z}_{\ge0}$.

(4) At $u=\infty$ the indicial equation fixes $\sigma=d+gm$, where $m=\deg h$. The exponent of $y$ in the term of
highest degree in $z$ is $d-\omega m\ge0$, which gives $d(1+\omega/g)\ge\omega\sigma/g$; since $1+\omega/g<0$
forces $g<0<\omega+g$, this is the stated bound.

(6) If some $r_i$ is a positive constant, if $1+\omega/g<0$, or if some $r_i$ is not a rational constant, then
(3), (4) or (5) leaves finitely many $d$. Otherwise every $r_i\in\mathbb{Q}_{\le0}$ and $1+\omega/g\ge0$, which is
exactly the condition for $H$ to be a polynomial. A homogeneous first integral $y^dh$ satisfies $\theta_dh=0$, so $h=\prod(u-u_i)^{-dr_i}$.
That $H$ is a first integral follows from $\X_{\rm top}(\log y)=y^eu=\X_{\rm top}\bigl(\int u\,du/G\bigr)$.
Conversely, the operator $\X_{\rm top}^2-\alpha\X_{\rm top}+\beta$ commutes with multiplication by first integrals
of $\X_{\rm top}$, so $H^kD^0_{\rm top}$ solves~(\ref{equ}) whenever $D^0_{\rm top}$ does.
\end{proof}

\begin{remark}[algorithmic scope of Theorem~\ref{thm:B}]\label{rem:Bscope}
Theorem~\ref{thm:B} is stated for the true $R$, and its bound uses $R_{[e]}$ through $\alpha$ and $\beta$. By
Proposition~\ref{prop:decouple}, $R_{[e]}$ is known up to the coefficient $c_e(x)$ of $y^e$ in the correction $c$,
and that coefficient depends on $D_{\rm top}$. Write $c_e=c^{\rm known}+g_1m+g\,s_1(h)$, where $g_1$ is the
coefficient of $u$ in $G$, $s_1$ is the sum of the roots of $h$ and $m=\deg h$. If $c_e$ is left free, the top level no longer bounds $d$. For H1 the top relation
becomes $d=1+3(n_1+n_2)$ with $c_e=-n_1-3n_2$, so every $(n_1,n_2)\in\mathbb{N}^2$ is compatible with the top.
Hence the bound of Theorem~\ref{thm:B} is a \emph{proven} bound for solutions with $c=0$, and the correction step
of the algorithm escalates heuristically; bounding $c_e$ from the next level of the expansion is open. For Theorem~\ref{thm:C} the correction has weight $0$, so it is a function $c_0(x)$ of $x$ alone. It shifts $K$
by $c_0(x)$, which changes neither $[K]_u$ nor the residues $k_j$ (because $f_jc_0/Y(f_j)\equiv0\pmod{f_j}$). Hence
the \emph{degree} statement of Theorem~\ref{thm:C} is unconditional. However, $c_0$ need not be a logarithmic
derivative. The $D$-step stays linear nevertheless: with $E=\exp(-\int c_0)$ one has
$c_0-\X(D)/D=-\X(ED)/(ED)$, so $S=R-\X(ED)/(ED)$ and $ED$ is a solution of the \emph{same} operator
$\mathcal{L}_R$, hyperexponential in $x$. In coordinates this is an exponential (hyperexponential) solution of the
linear differential system in $x$ of Section~\ref{sec:algo}, computable by the algorithm of~\cite{Pflugel1997}.
When $K$ has no poles, $c_0=h_0'/h_0$ is a logarithmic derivative and can be absorbed into $D$; by
Corollary~\ref{cor:Cfull} this is always the case when $N_0=1$, so the hyperexponential variant is needed only
for $N_0\neq1$.
\end{remark}

\begin{remark}[contrast with Kovacic]
In Kovacic's algorithm the degree of $P$ comes from a count of exponents (a Fuchs relation). Here that count is
empty: by the residue theorem $\sum r_i=1/g$, and $d$ drops out of it. The degree is bounded instead by the
\emph{positivity} of the multiplicities. With $R$ known, the only obstruction is a polynomial first integral of the
top field, which plays the role of resonance.
\end{remark}

\subsection{The regime $e=0$: an exact degree when the top is non-Liouvillian}
On the face $w=(0,1,1)$ all pieces of $\X$ have weight $0$. With $u=z/y$ and $\phi_{\rm top}=yF(x,u)$,
\[
\X_{\rm top}=\partial_x+u\,y\,\partial_y+G(x,u)\,\partial_u,\qquad G=F-u^2 ,
\]
so $u$ follows the \emph{Riccati equation} $u'=G(x,u)$ along $x$. Write $Y=\partial_x+G\partial_u$.

\begin{theorem}[Theorem C]\label{thm:C}
Assume that $F$ is polynomial in $u$, that $G$ has degree exactly $2$ in $u$ with leading coefficient $g_2\ne0$,
and that $u'=G(x,u)$ has no algebraic solution over $\overline{\K}(x)$; equivalently, the associated linear
equation has no Liouvillian solution, which is certified by Kovacic's algorithm. Then:
\begin{enumerate}
\item $S_{[0]}=-z/y$, and this is the unique rational solution of weight $\le0$ of the top Riccati equation;
\item with $K:=R_{[0]}+u$ (known up to an additive function of $x$, which affects neither $[K]_u$ nor the
$k_j$; see Remark~\ref{rem:Bscope}), $D_{\rm top}=y^d\,h_0(x)\prod_jf_j^{k_j}$, where the $f_j$ are the
irreducible factors in $u$ of the poles of $K$, $k_j\equiv f_jK/Y(f_j)\pmod{f_j}$, and
\[
d=[K]_u-g_2\sum_jk_j\deg_uf_j ,
\]
with $[K]_u$ the coefficient of $u$ in $K$ at $u=\infty$. In particular $d$ is unique; if some $k_j$ is not a
non-negative integer, no polynomial $D$ exists.
\end{enumerate}
\end{theorem}

\begin{proof}
(0) The weight of $S$ is at most $0$: with weight $\sigma>0$, level $2\sigma$ of~(\ref{ric}) would give
$S_{\rm top}^2=0$.

(1) In $u$, the top Riccati equation is $s_x+Gs_u=s^2+F_us-F+uF_u$, and $s=-u$ solves it for every $F$.
\begin{itemize}
\item $s=0$, which corresponds to weight $<0$, solves it only if $F=c(x)u$; then $u=0$ is a rational solution of
$u'=cu-u^2$, which is excluded.
\item Any other solution has the form $s=-u+1/v$, where $L(v):=v_x+(Gv)_u=-1$.
\begin{itemize}
\item A pole of order $k$ of $v$ along a non-invariant factor $f$ produces the uncancelled term
$-kP\,Y(f)/f^{k+1}$.
\item The invariant factors that are polynomial in $u$ are the minimal polynomials of algebraic solutions of
$u'=G$, and there are none.
\item Hence $v\in\K(x)[u]$. If $v=a_ku^k+\cdots$ with $k\ge0$, the coefficient of $u^{k+1}$ in $L(v)$ is
$(k+2)g_2a_k\ne0$, which is impossible.
\end{itemize}
\end{itemize}
Note that the difference of two solutions of $L(v)=-1$ would be a rational Jacobi multiplier of $Y$.

(2) At level $0$, $S_{[0]}=R_{[0]}-\X_{\rm top}(D_{\rm top})/D_{\rm top}$. With $D_{\rm top}=y^dh$ this gives
$Y(h)/h=K-du$. The poles of $Y(h)/h$ are exactly the non-invariant factors of $h$, and there are no invariant
ones. As $u\to\infty$, $Y(f)/f\sim g_2(\deg_uf)\,u$; comparing coefficients of $u$ gives the formula for $d$.
\end{proof}

\begin{corollary}\label{cor:C}
If every factor $f$ occurring in $R$ has a top that is monomial in $y$, then $K$ has no poles and
$d=1+\sum_f\kappa_f\,w(f)$.
\end{corollary}

When $N_0=1$ this is always the case, and much more follows.

\begin{corollary}[the regime of Theorem~\ref{thm:C} with $N_0=1$]\label{cor:Cfull}
In the setting of Theorem~\ref{thm:C}, let $S$ be a rational solution of~(\ref{ric}) with $N_0=1$ satisfying
(H$_{\log}$), with exponents $\kappa_f$. Then:
\begin{enumerate}
\item $T=\operatorname{pp}_w(S)=0$;
\item $K$ has no poles, $D_{\rm top}=y^dh_0(x)$ and $d=1+\sum_f\kappa_f\deg_yf$;
\item the correction of Proposition~\ref{prop:decouple} is $c=h_0'/h_0$, so after replacing $D$ by $D/h_0$ one has
\[
S=R-\frac{\X(D)}{D},\qquad R=\sum_f\kappa_f\Bigl(\frac{\X(f)}{f}-\frac{a_f'}{a_f}\Bigr),
\]
where $a_f(x)y^{\deg_yf}$ is the top of $f$. Thus $R$ is determined by the exponents $\kappa_f$ alone.
\end{enumerate}
\end{corollary}
\begin{proof}
(1) By Theorem~\ref{thm:C}, $S$ has weight $\le0$ and $S_{[0]}=-z/y$. Polynomial monomials have weight $\ge0$, and
$-z/y$ is not one, so $\operatorname{pp}_w(S)=0$. (2) Since $N_0=1$, every factor $f$ of $N$ is characteristic,
$f=f(x,y)$, with top $a_f(x)y^{n_f}$, $n_f=\deg_yf$. Hence $\X(f)/f$ has weight-zero part $a_f'/a_f+n_fu$ and
$\operatorname{pp}_w(\X(f)/f)=a_f'/a_f$. So $R_{[0]}=\bigl(\sum\kappa_fn_f\bigr)u+c_0(x)$ and
$K=\bigl(1+\sum\kappa_fn_f\bigr)u+c_0$ has no poles; Theorem~\ref{thm:C}(2) gives $h=h_0(x)$ and the value of $d$.
(3) $\X(D)/D$ has weight-zero part $h_0'/h_0+du$, whose polynomial part is $h_0'/h_0$; and
$\X(D/h_0)/(D/h_0)=\X(D)/D-h_0'/h_0$.
\end{proof}

In this regime the heuristic steps (L1), (L2) and (L4) of Section~\ref{sec:algo} disappear: the top is proven to be
$-z/y$, there is no Newton recursion ($T=0$), and there is no correction. For each choice of exponents $\kappa_f$
(rational roots of~(\ref{indicial}), or $0$), the linear step with the $R$ of Corollary~\ref{cor:Cfull}(3) and
$w\deg D=d$, followed by the rational $x$-step of Section~\ref{sec:algo}, returns \emph{all} such $D$. When no indicial equation is satisfied identically in $\kappa$, there are finitely
many choices, and the procedure is therefore a \emph{decision procedure} for rational $S$ with $N_0=1$ satisfying (H$_{\log}$); when $N$ has no
characteristic factor with a pole of $S$ (for instance $N=1$), (H$_{\log}$) is vacuous and it decides the existence
of rational $S$ outright. An empty kernel for every choice of $\kappa$ is a proof of non-existence
(Section~\ref{sec:negative}).

\begin{remark}
Neither hypothesis of Theorem~\ref{thm:C} can be dropped. (i) $\deg_uG=2$ is needed: for
$y''=z^2/y+y/x$ ($G$ of degree $0$) and $y''=z^2/y+z+y/x$ (degree $1$), the top Riccati equation has no algebraic
solution, yet there is a second rational $S$ of weight $0$. (ii) With a Liouvillian top the conclusion fails:
with the invariant factor $z-xy$, solutions with $D_{\rm top}=y(z-xy)^n$ exist and give $d=n+1$ for $n=1,2,3$,
while the formula gives $1$.
\end{remark}

\section{A certified search procedure}\label{sec:algo}

Algorithm~\ref{alg:KS} collects the steps of the previous sections. It is a search over finitely many local
choices, $\kappa_f$ and the top data at infinity, followed by linear (or reduced bilinear) algebra for $D$. Every
output is certified by substitution in~(\ref{ric}), so the algorithm has no false positives.

\paragraph{The $x$-dependence of $D$.} Theorems~\ref{thm:B} and~\ref{thm:C} use weights with $w_x=0$, so they bound
the degree of $D$ in $(y,z)$ only. The dependence on $x$ needs no bound. Write $D=\sum_jd_j(x)\,m_j(y,z)$ over the
finitely many monomials of weight at most $d_{\max}$. Then~(\ref{LR}) becomes a linear differential system in $x$
for $(d_j)$, and its rational solutions are computed, with no a-priori bound on their degree, by the algorithms of
Abramov and Barkatou~\cite{Abramov1989,Barkatou1999}. These systems are overdetermined and of mixed order (in H4, 32
equations for 10 unknowns). We reduce them to first order ($e_j=d_j'$), solve the purely algebraic equations, and
apply \texttt{RationalSolution} (Maple's \texttt{LinearFunctionalSystems}). On H4 and H5 this routine stops with an
internal error; there we first bring the system to normal form by linear differential elimination (prolongation and
Gaussian elimination over $\K(x)$, as in \texttt{rifsimp}), after which the unknowns satisfy a small normal system.
Table~\ref{tab:xdeg} shows the outcome: in all nine cases of battery~H the general rational solution is $\K\cdot D$.
For H4 and H5 the normal form is small enough to read off by hand. For
H5 it is $d_2'=0$, $d_1''=xd_1$, $d_0=d_1+xd_2$ (in the basis $1,y,y^2,\dots$); Airy's equation has no nonzero rational
solution, so $D=c\,(x+y^2)$. For H4, eliminating $d_3$ gives $w'''=4xw'+2w$ for $w=d_6$, the symmetric square of
Airy's equation; it has no finite singularity and a polynomial solution of degree $n$ would need $4n+2=0$, so $w=0$,
then $d_3$ is constant and $d_0-d_3$ solves Airy's equation, whence $D=c\,(y^3+xy+1)$. (Maple's \texttt{ratsols}
confirms that neither Airy's equation nor its symmetric square has a nonzero rational solution.) The non-Liouvillian top of
Theorem~\ref{thm:C} reappears here as the obstruction to further solutions.

\paragraph{Limits.} Completeness is limited by the following steps, which are the only heuristic ones.
\begin{enumerate}
\item[(L1)] The top equation at infinity is solved with an ansatz of quasi-homogeneous Laurent monomials, and
degenerate tops and resonances of the Newton recursion are not treated (Section~\ref{sec:infinity}). On balanced
faces with $e>0$ this costs speed, not solutions (Remark~\ref{rem:balanced}); it matters only on non-balanced
faces.
\item[(L2)] The Newton recursion for $T$ uses a bounded window in the $x$-degree of its coefficients.
\item[(L3)] The local step uses the indicial equation, which presupposes (H$_{\log}$).
\item[(L4)] For $e>0$ the weighted degree of $D$ is proven only for solutions whose correction $c$ vanishes
(Remark~\ref{rem:Bscope}); outside the regimes of Theorems~\ref{thm:B} and~\ref{thm:C}, and for $N_0\neq1$, the
degree limit is heuristic.
\end{enumerate}
A negative answer is therefore a proof of non-existence only for solutions satisfying the hypotheses behind
(L1)--(L4). In the regime of Theorem~\ref{thm:C} with $N_0=1$, only (L3) remains
(Corollary~\ref{cor:Cfull}).

\begin{algorithm}[ht]
\caption{Kovacic-$S$}\label{alg:KS}
\begin{algorithmic}[1]
\Require $\phi=M/N\in\K(x,y,z)$ in lowest terms.
\Ensure A list of families $S_\lambda$ of rational solutions of~(\ref{ric}), each certified, with a flag
(\emph{proven bound} or \emph{heuristic bound}). Every output solves~(\ref{ric}).
\State Factor $N=\prod f^{m_f}$; classify each $f$ as characteristic ($f_z=0$) or not; set
$N_0=\prod_{f_z\ne0}f^{m_f}$.
\State Under (H$_{\log}$): for each characteristic $f$, compute the rational roots of the indicial
equation~(\ref{indicial}); add $\kappa_f=0$.
\State Enumerate the faces of the Newton polytope of $(M,N)$, balanced faces first. For each face and each admissible
$\sigma$ in~(\ref{balance}), solve the top equation and run the linear recursion. Collect candidate pairs
$(S_{\rm top},T)$.
\For{each choice of $(\kappa_f)$ and $(S_{\rm top},T)$}
  \State $R\gets T-\operatorname{pp}_w\bigl(\sum\kappa_f\X(f)/f\bigr)+\sum\kappa_f\X(f)/f$ \Comment{linear step: $c=0$}
  \If{$N_0=1$}
     \State Compute $d_{\max}$ from Theorem~\ref{thm:B} ($e>0$) or Theorem~\ref{thm:C} ($e=0$, after certifying
     the top with Kovacic's algorithm); if neither applies, set a heuristic limit.
     \State Linear step: write $D=\sum_jd_j(x)m_j$ over the monomials with $w\deg\le d_{\max}$; reduce
     $\mathcal{L}_R(D)=0$ to normal form by linear differential elimination and compute its general rational
     solution (\texttt{RationalSolution}); the whole kernel is returned, with no early stop. For $e=0$ with poles
     in $K$ (possible only if $N_0\neq1$), compute instead the hyperexponential solutions $E(x)D$
     (Remark~\ref{rem:Bscope}).
  \Else
     \State Linear step as above (with $P_0=0$).
     \State Bilinear step: the reduced bilinear system in $(P_0,D)$, $\deg P_0<\deg N_0$, by Gr\"obner bases,
     raising $\deg D$ from $1$.
  \EndIf
  \State Correction step, if nothing was found: the correction $c$ of weight in $[0,e]$ unknown
  (Proposition~\ref{prop:decouple}), which is bilinear in $(c,D)$, with a heuristic degree limit. Run from $T=0$
  it covers every polynomial part on balanced faces (Remark~\ref{rem:balanced}).
  \State Structural step (Theorem~\ref{thm:struct}, Remark~\ref{rem:N0}): general $P/N$ over \emph{all} factors of
  $N$, with $\operatorname{pp}_w(P/N)=0$, plus the options $-\X(f)/f$ on non-characteristic factors and $D$. This
  step is bilinear, with a heuristic degree limit.
  \State For each kernel, form $S_\lambda=R-\X(D_\lambda)/D_\lambda$ with $D_\lambda=\sum\lambda_iD_i$, and certify
  it in~(\ref{ric}) with symbolic $\lambda$.
\EndFor
\State Flag the trivial branch $S=-\phi/z$ when $\phi_x=0$, and return the families.
\end{algorithmic}
\end{algorithm}

\begin{proposition}[soundness and relative completeness]\label{prop:alg}
Every family returned by Algorithm~\ref{alg:KS} consists of solutions of~(\ref{ric}). Conversely, let $S$ be a
rational solution with $N_0=1$, satisfying (H$_{\log}$), whose top $S_{\rm top}$ and polynomial part $T$ are among
the candidates of step~3. Assume that either the regime of Theorem~\ref{thm:C} applies, or the regime of
Theorem~\ref{thm:B} applies with one of its items (3)--(5) giving a finite $d_{\max}$ and $c=0$
(Proposition~\ref{prop:decouple}). Then $S$ belongs to one of the returned families, and the flag is \emph{proven
bound}.
\end{proposition}
\begin{proof}
Soundness holds by certification. For completeness: by Theorem~\ref{thm:den} and (H$_{\log}$),
$S=R-\X(D)/D$ with $R$ equal to the $R$ of the linear step for the correct choice of $\kappa_f$ and $T$ (up to
the additive $c_0(x)$ when $e=0$, which is absorbed by the hyperexponential factor $E$). Then $D$ (or $ED$) lies in
the kernel of $\mathcal{L}_R$, its weighted degree is at most $d_{\max}$ by Theorem~\ref{thm:B} or~\ref{thm:C}, and
its coefficients are rational (or hyperexponential) solutions of the system in $x$, which the rational
(exponential) solution algorithm returns completely.
\end{proof}
The hypotheses ``top among the candidates'' and ``$c=0$'' are where the procedure is heuristic ((L1)--(L4)); they
are what separates it from a decision procedure.

\paragraph{Design choices that mattered.}
\begin{enumerate}
\item \emph{The output is a kernel, not a single $S$.} The map $D\mapsto-\X(D)/D$ is not linear, so a
kernel of dimension $\ge2$ yields a family that is not linear in the parameters. For $y''=-y$ the basis
$\{z,y\}$ gives $S_\lambda=(\lambda_1y-\lambda_2z)/(\lambda_1z+\lambda_2y)$, a projective line of solutions, of
which $S=(y-z)/(y+z)$ is the member $\lambda=(1,1)$.
\item \emph{Linear before bilinear.} All candidates $(\kappa,T)$ are first tried with the cheap linear step, and
only then with the bilinear ones. Raising $\deg D$ before exhausting the candidates, or interleaving the two steps for each
candidate, made the budget evaporate in our tests (H6 and the blind cases C5, C6 failed).
\item \emph{Spurious factors.} Before computing a Gr\"obner basis, factors of the numerator of the equations that do
not involve the unknowns are removed (66 equations become 28 in one example); otherwise the ideal explodes.
\item \emph{The engine.} The cost of the bilinear step is a matter of implementation, not of arithmetic.
\begin{itemize}
\item Computing over GF$(p)$ costs the same as over $\mathbb{Q}$ in sympy's Buchberger implementation (108.6\,s
against 107.9\,s, for 13 unknowns): coefficient growth is not the bottleneck there.
\item Maple's \texttt{Groebner[Basis]} solves the same systems two or three orders of magnitude faster
(Table~\ref{tab:bilinear}).
\end{itemize}
\end{enumerate}

\section{Experiments}\label{sec:exp}

All computations were done twice, in sympy (Python) and in Maple~17, by independent code. The tables report the
Maple times on an AMD Ryzen~7 5000-series processor under Windows. Scripts and outputs are in the supplementary
material.

\subsection{Battery H: non-Liouvillian pull-backs}
The nine equations of battery~H are pull-backs of Airy's and Bessel's equations by Riccati substitutions
$B=z/y+c(x,y)$ and point transformations. All of them have non-Liouvillian first integrals~\cite{PaperA}. For H4 the
absence of regular $S$ follows from Corollary~\ref{cor:nonhyper}; for the others we only observe that the regular
ansatz finds none. In all of them $N_0=1$ and $c=0$, and (H$_{\log}$) holds: every characteristic factor along which
$S$ has a pole is simple, the polar part is logarithmic with constant $\kappa$, and $D$ is free of $z$ (in H6, $S$ has
no pole on the double factor $x$). So only the linear step is used.

\begin{table}[ht]\centering\footnotesize\setlength{\tabcolsep}{4pt}
\caption{Battery H. The columns give $N$, the polynomial $D$ found, the exponents $\kappa$, the weighted degree
bound and its source (bounds from Theorem~\ref{thm:B} are conditional on $c=0$), and the time of the linear step
(Maple). The sympy prototype took 0.8--141\,s per case.}\label{tab:H}
\begin{tabular}{@{}llllll@{}}\toprule
case & $N$ & $D$ & $\kappa$ & $w\deg D$ (bound) & $D$-step\\\midrule
H1: $c=y^2$, Airy & $1$ & $y$ & --- & $1$ (Thm.~B, $\omega=3$) & 0.0\,s\\
H2: H1$\circ(y^3+xy+1)$ & $3y^2+x$ & $y^3+xy+1$ & $1$ & $3$ (Thm.~B, $\omega=7$) & 0.3\,s\\
H3: H1$\circ(y+y^3)$ & $3y^2+1$ & $y(y^2+1)$ & $1$ & $3$ (Thm.~B, $\omega=7$) & 0.1\,s\\
H4: Airy$\circ(y^3+xy+1)$ & $3y^2+x$ & $y^3+xy+1$ & $1$ & $=3$ (Thm.~C) & 0.0\,s\\
H5: Airy$\circ\frac{y^2+x}{y+1}$ & $(y+1)(x-y^2-2y)$ & $x+y^2$ & $-1,\,1$ & $=2$ (Thm.~C) & 0.7\,s\\
H6: $c=y^2$, Bessel$_0$ & $x^2$ & $y$ & $0$ & $1$ (Thm.~B, $\omega=3$) & 0.0\,s\\
H7: $c=y^3$, Airy & $1$ & $y$ & --- & $1$ (Thm.~B, $\omega=4$) & 0.0\,s\\
H8: H7$\circ(y+xy^2)$ & $2xy+1$ & $y(xy+1)$ & $1$ & $2$ (Thm.~B, $\omega=7$) & 0.2\,s\\
H9: $c=xy^2$, Airy & $1$ & $y$ & --- & $1$ (Thm.~B, $\omega=3$) & 0.0\,s\\\bottomrule
\end{tabular}
\end{table}

\begin{itemize}
\item The $S$ of all nine are found and certified, with the same $\kappa$, $T$, $D$ and $S$ in both implementations.
\item All nine degree bounds are attained. Those from Theorem~\ref{thm:B} are conditional on $c=0$, which holds
here a posteriori.
\item For H4 and H5, Maple's \texttt{kovacicsols} returns no Liouvillian solution for the top linear equation
$w''=xw$, which certifies the hypothesis of Theorem~\ref{thm:C}. In both, $T=0$, $c=0$ and $R$ is the one of
Corollary~\ref{cor:Cfull}(3), as predicted. For H4 the indicial equation on $x+3y^2$ has the single root $\kappa=1$;
the choice $\kappa=0$ ($d=1$) gives the zero kernel, and $\kappa=1$ ($d=3$) gives $\K\cdot(y^3+xy+1)$
(Table~\ref{tab:xdeg}). Hence the $S$ of H4 is the \emph{only} rational $S$ with $N_0=1$ satisfying (H$_{\log}$).
\item A Gr\"obner computation on an ansatz without the structure of Theorem~\ref{thm:struct} (numerator and $D$
unknown together) did not finish H2 in 30 minutes.
\end{itemize}
Controls:
\begin{itemize}
\item $y''=z^2/y$ gives $S=-z/y$ ($D=1$, $\kappa=-1$);
\item Painlev\'e~I gives nothing;
\item the regular Liouvillian examples of~\cite{ACDM2019} (such as $S=-z/(x^4+y)$) are \emph{not} found, as
expected: they are regular, and the regular ansatz is the right tool for them.
\end{itemize}

\subsection{Coverage on a second non-Liouvillian family}\label{sec:coverage}
To see how often the linear step suffices outside battery~H, we generated 14 equations with the generator of
\cite{PaperA}. For each, $\phi=\bigl[x-g^2-g_x-zg_y\bigr]/g_z$ and $S=g_y/g_z$, with Airy's equation behind and
hence a non-Liouvillian first integral, for simple choices of $g$ (one more choice was discarded because the
classification script failed on it). All 14 have the form~(\ref{struct}).
\begin{center}\footnotesize
\begin{tabular}{@{}p{0.3\linewidth}p{0.36\linewidth}p{0.28\linewidth}@{}}\toprule
class & cases & how found\\\midrule
(H$_{\log}$), $N_0=1$ & 7 ($g=z/y+y^2$, $z/y+xy$, $yz+x$, \dots) & linear step\\
regular, with a non-logarithmic polar part & 4 ($g=z/(y^2+x)$, $(z+y^2)/(y+1)$, \dots) & regular ansatz / structural step with $D=1$\\
singular, needing $P/N$ and $D$ & 2 ($g=z/y+1/(x+y)$, $(z+xy^2)/(y(y+1))$) & structural step (70\,s and 18\,s, sympy)\\
$N_0\neq1$ & 1 ($g=z^2/y$) & bilinear step\\\bottomrule
\end{tabular}
\end{center}
So in this family half of the solutions satisfy the hypotheses of the linear step, and all of them are reached by
the structural step. In Maple the structural step finds the two singular cases in 0.12\,s and 0.17\,s (against 18\,s and 70\,s
in sympy), each with a single component that reproduces the true $S$.

\begin{table}[ht]\centering\footnotesize
\caption{The $x$-dependence of $D$: the general rational solution in $x$ of $\mathcal{L}_R(D)=0$ over the monomials
of weight $\le d_{\max}$, with $R$ the $R$ of the linear step. Maple~17 (\texttt{RationalSolution}); for H4 and H5,
linear differential elimination first (Maple \texttt{rifsimp}, which returns the same normal form as our sympy
elimination; times include it; sympy: 5.7\,s and 29\,s).}\label{tab:xdeg}
\begin{tabular}{@{}llll@{}}\toprule
case & $d_{\max}$ & general rational solution & time\\\midrule
H1, H6, H7, H9 & 1 & $c\,y$ & 0.02--0.14\,s\\
H2 & 3 & $c\,(y^3+xy+1)$ & 0.03\,s\\
H3 & 3 & $c\,y(y^2+1)$ & 0.22\,s\\
H8 & 2 & $c\,y(xy+1)$ & 0.30\,s\\
H4 & 3 & $c\,(y^3+xy+1)$; normal system of dimension 6 & 0.11\,s\\
H5 & 2 & $c\,(x+y^2)$; normal system of dimension 3 & 0.02\,s\\\midrule
Example~\ref{ex:gap}, $R=0$ & 3 & $c_2(xz+y/2)+c_1z$ (contains the trivial branch) & 0.09\,s\\\bottomrule
\end{tabular}
\par\smallskip\raggedright\footnotesize The last row uses $R=0$, the $R$ produced by the correction step, not the
$R=T$ of the linear step (for which the kernel is zero); it only illustrates the $x$-step.
\end{table}

\begin{table}[ht]\centering\footnotesize
\caption{Coverage of the search procedure.}\label{tab:coverage}
\begin{tabular}{@{}p{0.42\linewidth}p{0.5\linewidth}@{}}\toprule
type of solution & status\\\midrule
$N_0=1$, (H$_{\log}$), top among the candidates, and Theorem~\ref{thm:C}, or Theorem~\ref{thm:B} with $c=0$ &
found, with a proven weighted degree bound and the complete $x$-dependence (Prop.~\ref{prop:alg})\\
$c\neq0$ (Example~\ref{ex:gap}) & found by the correction step; heuristic bound\\
(H$_{\log}$) fails: higher-order or non-logarithmic poles, or $D$ meeting a characteristic factor
(Remark~\ref{rem:Hlog}) & form guaranteed by Theorem~\ref{thm:struct}; found by the structural step; heuristic bound\\
$N_0\ne1$ & bilinear step; heuristic bound\\
top not among the Laurent-monomial candidates, balanced face & found by the correction step from $T=0$
(Remark~\ref{rem:balanced}); heuristic bound\\
top not among the candidates on a non-balanced face, resonant recursion & may be missed ((L1), (L2))\\\bottomrule
\end{tabular}
\end{table}

\subsection{Certified negative answers}\label{sec:negative}
Corollary~\ref{cor:Cfull} turns an empty output into a proof. Two examples, both with Airy's equation as top
(certified by \texttt{kovacicsols}) and $T=0$:
\begin{itemize}
\item \textbf{N1}: $y''=xy+1$. Here $N=1$, so (H$_{\log}$) is vacuous, $R=0$ and $d=1$: $D=d_0+d_1y+d_2z$. The
system in $x$ reduces to $d_2''=-2d_1'$, $d_1''=-d_2-2xd_2'$, $d_0''=xd_0-d_1-2d_2'$. Eliminating $d_1$ gives the
symmetric square of Airy's equation for $d_2$, so $d_2=0$; then $d_1=a$ is constant and $d_0''-xd_0=-a$, whose
rational (hence polynomial) solutions must vanish by degree, forcing $a=0$. So the kernel is zero and
\emph{$y''=xy+1$ has no rational $S$}, although it has a (non-Liouvillian) first integral.
\item \textbf{N2}: the pull-back of N1 by $y=Y^3+xY+1$,
\[
y''=\frac{x^2y+xy^3+x-6yz^2-2z+1}{x+3y^2},
\]
which differs from H4 only by the term $1/(x+3y^2)$. The indicial equation on $x+3y^2$ again has the single root
$\kappa=1$. For $\kappa=0$ ($d=1$) the kernel is zero. For $\kappa=1$ ($d=3$, $R=\X(f)/f$) the normal form is that of
H4 except for $d_0''=xd_0-(x+1)d_3-\tfrac13d_6-\tfrac23(x+1)d_6'$; as for H4, $d_6=0$ and $d_3=a$ is constant, and
then $u=d_0-a$ satisfies $u''=xu-a$, so $u=0$ and $a=0$. Hence N2 has no rational $S$ with $N_0=1$ satisfying
(H$_{\log}$), while its neighbour H4 has exactly one.
\end{itemize}
The pull-back by a point transformation does not preserve rationality of $S$ in either direction (compare
Theorem~\ref{thm:A}), so N2 is not a consequence of N1. Both computations were repeated in Maple~17, where
\texttt{rifsimp} returns the same normal forms and \texttt{RationalSolution} the zero kernel, in 0.02--0.03\,s per
choice of $\kappa$ (with H4, $\kappa=1$, as positive control: $\K\cdot(y^3+xy+1)$ in 0.11\,s).

\subsection{Robustness tests (Liouvillian, $N_0\ne1$)}\label{sec:robust}
The bilinear $(P_0,D)$ step acts only when $N_0\neq1$, and the correction step only when $c\neq0$; neither happens in
battery~H. We tested them on
Liouvillian families built specifically for this purpose. They are outside the scope of the theorems, and serve only
to test the machinery.
\begin{itemize}
\item \textbf{Blind battery} C0--C9 (10 singular cases from random closed 1-forms; C5 is~(\ref{C5})): 10/10,
with 0--9\,s per case. The two
cases whose only solution is the trivial branch take about 160\,s; that is the cost of a negative answer.
\item \textbf{$N_0\ne1$ battery} (16 cases): 16/16, about 8\,s each (sympy).
\item \textbf{Between forms} (four point transforms, among them X0 with $T=Y^2$ and X2 with $T=xY$, raising
$\deg N_0$ to 4 and $\deg D$ to 2): 4/4.
\end{itemize}
In all three there were no false positives.

\begin{table}[ht]\centering\footnotesize
\caption{The reduced bilinear step: sympy (Buchberger) against Maple~17 (\texttt{Groebner[Basis]}). C5 is~(\ref{C5});
X0, X2 are from the ``between forms'' battery.}\label{tab:bilinear}
\begin{tabular}{@{}lccc@{}}\toprule
system & unknowns & sympy & Maple\\\midrule
C5, $\deg D=1$ & 6 & 2.3\,s & 0.05\,s\\
C5, $\deg D=2$ & 9 & 25\,s & 0.08\,s\\
C5, $\deg D=3$ & 13 & $>212$\,s (failed) & 0.30\,s\\
X0 ($T=Y^2$), $\deg D=1$ & $10+3$ & failed (400\,s) & 0.07\,s\\
X2 ($T=xY$), $\deg D=2$ & $20+10$ & --- & 8.0\,s\\\bottomrule
\end{tabular}
\end{table}

\subsection{Checks of the theory}\label{sec:checks}
\begin{itemize}
\item \emph{Theorem~\ref{thm:A}}: 41 sympy checks and Maple~17; a blind Gr\"obner search with an ansatz
$S=P/T_Y$ independent of the theorem (8/8 as predicted, in both systems); two larger blind searches with
separately written code, in an AI-assisted review session (24 and 32 unknowns; see the section on AI tools), which found exactly the predicted $S$, and new cases of branch~(ii) with
$c\neq0$ (for example $y''=2y/x^2$ with $T=Y^2-x^2$).
\item \emph{Local analysis}: Proposition~\ref{prop:char} on 53 factors of the group's examples and their transforms,
with a negative control (perturbing $S$ makes (E1) fail), and on 21 further checks with $m\ge2$;
Proposition~\ref{prop:nonchar} with generic functions $Z,\mu,s_j$ for $m=1,2$ up to $j=3$, and statement~(a) on 40
random factors.
\item \emph{Theorems~\ref{thm:B} and~\ref{thm:C}}: on battery~H in sympy and Maple~17, every bound attained. Two
separate AI-assisted review sessions, with their own code (see the section on AI tools), added about 70 exact top-level families with $n_i\in\{1,2\}$ and 12 complete examples,
including irrational and complex residues and a resonant case, plus 8 cases with a non-Liouvillian top. No
counterexample was found.
\end{itemize}
A pole in $K$ (Theorem~\ref{thm:C}(2) with $k_j\neq0$) cannot occur when $N_0=1$ (Corollary~\ref{cor:Cfull}); it
requires non-characteristic factors, which are outside the scope of the degree results.

\section{Related work and limitations}\label{sec:limits}

\paragraph{Kovacic-type algorithms.} Kovacic's algorithm \cite{Kovacic1986} decides Liouvillian solvability of
second-order linear ODEs. Its structure (local data at the singularities and at infinity, then a polynomial found
by linear algebra) is the template of the present work. Refinements of the local-to-global step include
\cite{vHW2005} and, for partial linear difference equations, the partial denominator bounds of \cite{KS2010}.
Rational and exponential solutions of linear differential systems, which we use for the $x$-dependence of $D$, are
computed by the algorithms of \cite{Abramov1989,Barkatou1999,Pflugel1997}. The scalar fact behind our residue lemma (Lemma~\ref{lem:res}(a),(c)) is classical: it appears in Kovacic's
analysis of case~1 \cite[Sec.~3]{Kovacic1986}; see also \cite{Echi}.

\paragraph{Darboux, Prelle--Singer and Riccati first integrals.} Prelle and Singer \cite{PS1983} and Singer
\cite{Singer1992} characterise elementary and Liouvillian first integrals of planar fields. Man and MacCallum
\cite{MM1997} made the planar algorithm effective. For algebraic solutions and invariant curves see
\cite{GGG2006}; for degree bounds on Darboux polynomials see Ch\`eze \cite{Cheze2011}. Ch\`eze and Combot \cite{CC2020} compute first integrals of polynomial vector fields, including Riccati-type
ones, by linear algebra with a degree bound as \emph{input}. They note that computing such a bound ``is still an open problem'' (Poincar\'e's problem). Our Theorems~B and~C give bounds of this kind for the
global factor $D$ of $S$, in the non-Liouvillian regime and under explicit hypotheses. In the language of
foliations, a rational $S$ makes the codimension-one foliation $\omega=-(\phi+zS)dx+S\,dy+dz$ integrable; pull-backs
of Riccati foliations are the transversely projective ones \cite{CLNPT2007,Casale2006}.

\paragraph{The $S$-function method and its neighbours.} The $S$-function method was introduced in
\cite{ACDM2019} and extended in \cite{NewS2023}. The same function $S$ appears in the extended Prelle--Singer method
of Chandrasekar, Senthilvelan and Lakshmanan \cite{CSL2005}, and $S$ is related to $\lambda$-symmetries \cite{MR2012}.
In \cite{CSL2005} the pair $(S,R)$, with $R$ an integrating factor, is determined from the Riccati equation for
$S$ and a linear equation for $R$. Both are found by ans\"atze chosen case by case, for $S$ a M\"obius form
$S=(a+bz)/(c+dz)$ with coefficients depending on $(x,y)$. Such ans\"atze admit denominators outside $N$, but the
location and order of the poles of $S$ and bounds on its degree are not discussed there; these are what the present
paper supplies. Huang
\cite{Huang2025,Huang2026} reduces rational ODEs to Riccati equations by different means. Gin\'e and
Sinelshchikov \cite{GS2025} construct non-Liouvillian integrals for autonomous equations by the same mechanism
that underlies \cite{PaperA}, and Demina and Nechitailo \cite{DN2026} study the $R$-integrability of such equations.

\paragraph{Limitations.}
\begin{itemize}
\item The heuristic steps are (L1)--(L4) of Section~\ref{sec:algo}: the Laurent-monomial ansatz for the top and
the unresolved degenerate or resonant cases of the Newton recursion; the $x$-window for $T$; (H$_{\log}$); and, for
$e>0$, the vanishing of the correction $c$ (bounding $c_e$ from the next level of the expansion is open). The degree
results assume $N_0=1$.
\item No a-priori bound on the $x$-degree of $D$ is proved, and none is needed: the $x$-dependence is computed by
rational-solution algorithms for linear systems. With Maple~17 this required a preliminary differential elimination
on H4 and H5.
\item Theorem~\ref{thm:B} needs $\deg F\le2$, $g\ne0$, simple roots of $G$ and a monomial $N_{\rm top}$.
Theorem~\ref{thm:C} needs $\deg_uG=2$ and a top Riccati equation without algebraic solutions; neither can be
dropped.
\item Not covered by the degree results: a Liouvillian top in the regime $e=0$ and the resonant top (a polynomial
first integral of $\X_{\rm top}$). In both a Darboux structure is present at the top level, and the regular ansatz of
\cite{ACDM2019,NewS2023} or Prelle--Singer-type methods are the natural tools. Factors of $N$ depending on $x$ alone
have the local analysis of Remark~\ref{rem:xonly}, which the implementation does not yet use.
\item The linear step excludes the options listed in Remark~\ref{rem:Hlog}; the structural step covers them, but
without a degree bound.
\item Theorem~A covers point transforms with $x$ fixed; general Lie point transformations, which also change $x$,
are not treated.
\end{itemize}

\section{Conclusion}
The rational $S$-functions that the regular ansatz cannot find are not an anomaly. For point transforms of linear
equations without hyperexponential solutions they are the generic case (Theorem~A), and in general they have a
precise structure. All of the
denominator is read off from $N$, except for one squarefree polynomial $D$ whose zeros are poles with residue
$-1$. The polynomial part comes from the Newton polygon at infinity, up to a correction of weight at most $e$. The
weighted degree of $D$ is determined exactly when the top is a non-Liouvillian Riccati equation, and bounded by
residues of the top field otherwise, conditionally on that correction. The resulting search transposes Kovacic's
strategy to the Riccati equation of $S$. It is certified and complete under stated hypotheses, but not yet a
decision procedure in general (it is one in the regime of Theorem~C with $N_0=1$, for solutions satisfying
(H$_{\log}$), and it proves, for instance, that $y''=xy+1$ has no rational $S$): the top ansatz and the resonant
recursion at infinity, the coupled local condition when
(H$_{\log}$) fails, and a bound on the correction $c_e$ remain open. The dependence of $D$ on $x$, on the other hand,
is obtained completely from rational solutions of a linear system in $x$; in the Airy cases it is the symmetric square
of the top equation that excludes further solutions. Three phenomena have no scalar counterpart: free Cauchy data on non-characteristic surfaces, the vanishing
of the Fuchs count by the residue theorem, and a top that is itself a Riccati equation. Together with \cite{PaperA}, which turns a rational $S$ into a decision about Liouvillian integrability and
recovers the underlying linear equation, this gives a working pipeline for the non-Liouvillian pull-backs considered
here.

\section*{Supplementary material}

The supplementary archive contains both implementations used in Section~\ref{sec:exp}: the \texttt{sympy}
prototype of the search (with the drivers and inputs of battery~H, of the blind and $N_0\neq1$ batteries and of
the coverage test) and the Maple~17 worksheets with their saved outputs (battery~H, the bilinear step, the
$x$-step of Table~\ref{tab:xdeg}, the certified negatives, the correction step, the balanced-face example and the
checks of Theorems~A, B and~C). A README file maps each script to the statement or table it supports.

\section*{Use of AI tools in the research}

AI assistants were used in this work, and we document their role here. Claude (Anthropic) was used in working
sessions to propose and check constructions and proofs, to write the \texttt{sympy} verification scripts and the
Maple worksheets, and to act as an internal referee of successive drafts. Separate review sessions of this kind,
with code written independently of the main implementation, produced the additional checks reported in
Section~\ref{sec:checks}. Every statement and proof was checked by the authors, and every computation reported in
Section~\ref{sec:exp} was reproduced by two independent implementations (Maple and \texttt{sympy}); the scripts are
in the supplementary material. The AI systems are not authors.

\section*{Declaration of generative AI and AI-assisted technologies in the manuscript preparation process}

During the preparation of this work the authors used Claude (Anthropic) in order to draft parts of the text and
to improve its language. After using this tool, the authors reviewed and edited the content as needed and take
full responsibility for the content of the published article.

\section*{Acknowledgements}

We thank Claude (Anthropic), the AI assistant described above, for its help in developing this work.

\section*{Funding}

\TODO{CNPq / FAPERJ: agencies and grant numbers, as in the companion paper.}

\section*{Declaration of competing interest}

The authors declare that they have no known competing financial interests or personal relationships that could
have appeared to influence the work reported in this paper.

\section*{CRediT authorship contribution statement}

\textbf{L.\,G.\,S.\ Duarte} and \textbf{L.\,A.\,C.\,P.\ da Mota}: Conceptualization, Methodology, Formal analysis,
Software, Validation, Writing -- original draft, Writing -- review \& editing.

\section*{Data availability}

No datasets were generated or analysed. The Maple worksheets and their outputs, the \texttt{sympy} scripts, and
the machine-readable inputs of all experiments are provided as supplementary material.

\end{document}